\documentstyle[12pt,amstex]{article}

\newcommand{\lrw}{\longrightarrow}
\newcommand{\Map}{\longmapsto}
\newcommand{\beqa}{\begin{eqnarray*}}
\newcommand{\eeqa}{\end{eqnarray*}} 
\newcommand{\la}{\mbox{$\langle$}} 
\newcommand{\ra}{\mbox{$\rangle$}}

\newcommand{\As}{\mbox{${\rm Aut}_{\mathfrak g}$}} 
 
\newcommand{\Outs}{\mbox{${\rm Aut}_{D({\mathfrak g})}$}} 
\newcommand{\Autg}{{\rm Aut}_{\mathfrak g}}

\newcommand{\ea}{E_\alpha}

\newcommand{\hs}{\hspace{.2in}} 
 
\newcommand{\Ad}{{\rm Ad}}

\def \fg{{\mathfrak g}} 
\def \fh{{\mathfrak h}} 
\def \fk{{\mathfrak k}} 
\def \fp{{\mathfrak p}} 
 
\def \fa{{\mathfrak a}} 
\def \fb{{\mathfrak b}} 
\def \fn{{\mathfrak n}}

\def \fz{{\mathfrak z}}

\def \fu{{\mathfrak u}} 
 
\def \la{{\langle}} 
\def \ra{{\rangle}} 
\def \R{{\mathbb R}} 
 
\def \C{{\mathbb C}} 
 
\def \i{{\sqrt{-1}}}

\newcommand{\qed}{\begin{flushright} {\bf Q.E.D.}\ \ \ \ \ 
                  \end{flushright} }

\def \piANt{{\pi_{AN, s}}}

\newtheorem{thm}{Theorem}[section] 
\newtheorem{lem}[thm]{Lemma} 
\newtheorem{prop}[thm]{Proposition} 
 
\newtheorem{cor}[thm]{Corollary} 
\newtheorem{rem}[thm]{Remark} 
 
\newtheorem{nota}[thm]{Notation} 
\newtheorem{dfn}[thm]{Definition}

\newtheorem{dfn-nota}[thm]{Definition-Notation} 
\newtheorem{dfn-prop}[thm]{Definition-Proposition} 
\newtheorem{dfn-lem}[thm]{Definition-Lemma}

\numberwithin{equation}{section} 
 
\begin{document}

\title{Thompson's conjecture 
for real semi-simple Lie groups} 
\author{Jiang-Hua Lu\thanks{Research partially supported by 
NSF grant DMS-0105195, by HHY Physical Sciences Fund and
by the New Staff Seeding Fund   
at HKU.} \hspace{.04in} and  \hspace{.04in} 
Sam Evens\thanks{Research partially supported by 
NSF grant DMS-0096110.}} 
 
%\address{Department of Mathematics, The University of Notre Dame 
%\newline Department of Mathematics, The University of Hong Kong, Pokfulam, Hong Kong} 
 
%\email{evens.1@nd.edu, jhlu@maths.hku.hk} 
 
\maketitle 
 
%\tableofcontents 
\begin{center}
{\it Dedicated to Professor Alan Weinstein for his 60th birthday}
\end{center}

\bigskip 
\begin{abstract} 
A proof of Thompson's conjecture for real semi-simple Lie groups has been given by 
Kapovich, Millson, and Leeb. In this note, we give another proof of the conjecture 
by using a theorem of Alekseev, Meinrenken, and Woodward from symplectic 
geometry. 
\end{abstract} 
 
\section{Introduction} 
\label{sec_intro}

Thompson's conjecture for the group $GL(n,\C)$, which relates eigenvalues 
of matrix sums and singular values of matrix products,  
was first proved by Klyachko in \cite{klya:thompson}.  
In \cite{amw:linear}, by applying a Moser argument to certain 
symplectic structures, Alekseev, Meinrenken, and Woodward gave a proof 
of Thompson's conjecture for all quasi-split real semi-simple Lie groups. In \cite{KML:1},  
Kapovich, Millson, and Leeb have, among other things, 
proved Thompson's conjecture 
for an arbitrary semi-simple Lie group $G_0$. In this note, we give a different  
 proof of Thompson's conjecture for arbitrary  
semi-simple real groups by extending the proof of Alekseev, Meinrenken, and Woodward 
for quasi-split groups. In fact, we prove a stronger result, Theorem \ref{thm_main-intro}, 
which implies Thompson's conjecture. 
Let $G_0=K_0A_0N_0$ be an Iwasawa decomposition of $G_0$ and let 
$\fg_0=\fk_0 + \fp_0$ be a compatible  
 Cartan decomposition of the Lie algebra of $\fg_0$. 
Theorem \ref{thm_main-intro} asserts that for each $l\ge 1$, 
there is a diffeomorphism 
 $L:(A_0N_0)^l \to ({\fp_0})^l$ which relates the addition on $\fp_0$ with  
the multiplication on $A_0N_0$ and intertwines naturally defined 
$K_0$-actions. 
When $G_0$  is quasi-split, Theorem 
\ref{thm_main-intro} follows from  
results in \cite{amw:linear}.   
The key step in our proof of Theorem \ref{thm_main-intro} for an arbitrary $G_0$ 
is to relate an  
arbitrary real semi-simple Lie algebra $\fg_0$ to a quasi-split real form in its 
complexification. 
 
In  Section \ref{sec_the-conjecture}, we state Theorem \ref{thm_main-intro} and show 
that it implies Thompson's conjecture. 
Inner classes of real forms and quasi-split real forms are reviewed in Section  
\ref{sec_inner}. The proof of Theorem \ref{thm_main-intro} is given in Section  
\ref{sec_the-proof}. 
Since the version of the 
Alekseev-Meinrenken-Woodward theorem we present in this paper 
is not explicitly stated in   
\cite{amw:linear}, we give an outline of its proof in the Section \ref{sec_amw}, the Appendix.

\bigskip 
\noindent 
{\bf Acknowledgments.} 
Although we do not explicitly use any results from \cite{philip:loci}, our paper 
is very much inspired by \cite{philip:loci}. We thank P. Foth for showing us a  
preliminary version of \cite{philip:loci} and for  
helpful discussions. We also thank 
 J. Millson and M. Kapovich for sending us 
the preprint \cite{KML:1} and R. Sjamaar for answering some questions. 
 
\section{Thompson's conjecture} 
\label{sec_the-conjecture}

Let $G$ be a complex connected reductive algebraic group with an anti-holomorphic  
involution 
$\tau$. Let $G_0$ be a subgroup of the fixed point set $G^\tau$ of $\tau$ 
containing the identity connected component. Then $G_0$ is a 
real reductive Lie group in the sense of \cite{wallach:book} pp. 42--45, 
 which implies 
that $G_0$ has Cartan and Iwasawa decompositions. Let 
$\fg_0$ be the Lie algebra 
of $G_0$.  Fix a Cartan decomposition $\fg_0 = \fk_0 + \fp_0$ of $\fg_0$, 
and let $G_0 = P_0 K_0$ be the corresponding Cartan decomposition of $G_0$.    
Let $\fa_0 \subset \fp_0$ be a maximal abelian subspace of $\fp_0$. Fix a choice  
$\Delta_{res}^{+}$ of  
positive roots in the restricted root system $\Delta_{res}$ 
for $(\fg_0, \fa_0)$, and let $\fn_0$ be the subspace of $\fg_0$ spanned by the   
root vectors 
for roots in $\Delta_{res}^{+}$. 
Then 
 $\fg_0 = \fk_0 + \fa_0 + \fn_0$ is an Iwasawa decomposition for $\fg_0$. 
 Let $A_0 = \exp(\fa_0)$ and $N_0=\exp(\fn_0)$. 
Then we have the Iwasawa decomposition $G_0 = K_0 A_0 N_0$ for $G_0$. 
 
Consider now the space $X_0:=G_0/K_0$ with the left $G_0$--action given by 
\begin{equation} 
\label{eq_left-K0-action} 
G_0 \times (G_0/K_0) \lrw G_0/K_0: \, \, (g_1, \, gK_0) \Map g_1gK_0, 
\hs g_1, \, g \in G_0. 
\end{equation} 
Thompson's conjecture is concerned  with $K_0$--orbits in $X_0$. 
Identify 
\[ 
\fp_0 \, \stackrel{\exp}{\cong} \, P_0 \, \cong \, G_0/K_0 = X_0 
\] 
via the Cartan decomposition $G_0 = P_0K_0$. The $K_0$--action on $X_0$ in (\ref{eq_left-K0-action}) 
becomes the adjoint action of $K_0$ on $\fp_0$.  
Orbits of $K_0$ 
in $\fp_0$ are called (real) flag manifolds. Let $\fa_{0}^{+}\subset \fa_0$ be the 
Weyl chamber determined by $\Delta_{res}^{+}$. It is well-known that  
every $K_0$--orbit in $\fp_0$ goes through a unique element
 $\lambda \in \fa_{0}^{+}$.  
 
On the other hand, we can also identify  
\[ 
A_0N_0 \, \cong \, G_0/K_0 = X_0 
\] 
via the Iwasawa decomposition $G_0 = A_0N_0K_0$. Then the $K_0$--action on $X_0$ 
becomes the following action of $K_0$ on $A_0N_0$: 
\begin{equation} 
\label{eq_K0-on-A0N0} 
k \cdot b := p(kb) = p(kbk^{-1}), \hs {\rm for} \hs k \in K_0, \, b\in A_0N_0, 
\end{equation} 
where $p: G_0 \to A_0N_0$ is the projection $b_1k_1 \mapsto b_1$ for $k_1 \in K$ and  
$b_1 \in A_0N_0$. Let $E_0: \fp_0 \to A_0N_0$ be the composition of the identifications: 
\begin{equation} 
\label{eq_E0} 
E_0: \, \, \fp_0 \stackrel{\exp}{\cong} P_0 \cong G_0/K_0 \cong A_0N_0. 
\end{equation} 
Then $E_0$ is $K_0$--equivariant, and $E_0(\fa_0) = A_0$.  
Thus every $K_0$--orbit in $A_0N_0$ goes through a unique point $a = \exp \lambda \in A_{0}^{+} 
:=\exp \fa_{0}^{+}$. 
 
Thompson's conjecture for $G_0$ is concerned with the sum of $K_0$--orbits in $\fp_0$  
and the product of $K_0$--orbits in $A_0N_0$. To further 
prepare for the statement of the conjecture, let $l \geq 1$ be an integer, and consider the two maps 
\beqa 
& & {\bf a}: \,  
 \fp_0 \times \fp_0 \times \cdots \times \fp_0 \lrw \fp_0: \, (\xi_1, \xi_2, \cdots, \xi_l) \Map  
\xi_1 + \xi_2 +  
 \cdots + \xi_l,\\ 
& & {\bf m}: \, A_0N_0 \times A_0N_0 \times \cdots \times A_0N_0  \lrw A_0N_0: \, \,  
 \,  (b_1, b_2, \cdots, b_l) \Map b_1b_2 
 \cdots b_l. 
\eeqa 
Clearly,  
${\bf a}$ is $K_0$--equivariant for the diagonal  
action of $K_0$ on $(\fp_0)^l$. On the other hand, define the  
{\it twisted diagonal action} ${{\mathcal T}}$ of $K_0$ on $(A_0N_0)^l$ by 
\begin{equation} 
\label{eq_mathcalT} 
k \Map {{\mathcal T}}_k:=  
\nu^{-1} \circ \delta_k \circ \nu: \,\, (A_0N_0)^l \lrw (A_0N_0)^l, 
\end{equation} 
where $\delta_k$ is the diagonal action of $k \in K_0$ on $(A_0N_0)^l$, and  
$\nu:  (A_0N_0)^l \lrw (A_0N_0)^l$ is the diffeomorphism given by  
\begin{equation} 
\label{eq_nu} 
\nu(b_1, b_2, \cdots, b_l) \Map 
(b_1, \, \, b_1b_2,  \, \, \cdots, \, \, b_1b_2 \cdots b_l). 
\end{equation} 
 See Remark \ref{rem_twisted} for motivation of 
the twisted diagonal action. 
Let $e$ be the identity element of 
$A_0N_0$ and identify $T_{e}(A_0N_0) \cong \fa_0 + \fn_0 \cong \fg_0 /\fk_0 \cong \fp_0$. 
 We will regard the map ${\bf a}$, respectively the diagonal $K_0$--action on $(\fp_0)^l$, 
as the linearization of the map ${\bf m}$, respectively the twisted  
diagonal $K_0$--action on $(A_0N_0)^l$, 
 at the point 
$(e, e, \cdots, e)$.

\begin{nota} 
\label{nota_orbits} 
{\rm For $\lambda \in \fa_0$, we will use ${{O}}_\lambda$ to denote the 
$K_0$--orbit in $\fp_0$ through $\lambda$. For $a \in A_0$, we will use 
${{ D}}_a$ to denote the $K_0$--orbit in $A_0N_0$ through the point $a$.  
If $\Lambda = (\lambda_1, \lambda_2, \cdots, \lambda_l) \in (\fa_0)^l$, we set 
$a_j = \exp(\lambda_j)$ for $1 \leq j \leq l$, and 
\[ 
O_{\Lambda} = {{ O}}_{\lambda_1} \times {{O}}_{\lambda_2} \times \cdots \times  
{{ O}}_{\lambda_l} \hs {\rm and} \hs  
D_{\Lambda} = {{D}}_{a_1} \times  
{{D}}_{a_2} \times \cdots \times  
{{ D}}_{a_l}. 
\] 
} 
\end{nota} 
 
In this paper, we will  prove the following theorem. 
 
\begin{thm} 
\label{thm_main-intro} 
For every integer $l \geq 1$, there is a $K_0$--equivariant diffeomorphism  
$L: (A_0N_0)^l \to (\fp_0)^l$ such that ${\bf m} = E_0 \circ {\bf a} \circ L$. Moreover,  
$L(D_\Lambda) = O_\Lambda$ for every $\Lambda \in (\fa_0)^l$. 
\end{thm}

\bigskip 
Theorem \ref{thm_main-intro} now readily implies the following Thompson's conjecture for $G_0$. 
 
\begin{cor} [Thompson's conjecture] 
\label{thm_Thompson} 
For each $\Lambda = (\lambda_1, \lambda_2, \cdots, \lambda_l) \in (\fa_0)^l$,  
the two spaces 
\[ 
({\bf m}^{-1}(e) \cap D_\Lambda)/K_0 =  
\{(b_1, b_2, \cdots, b_l) \in D_\Lambda: b_1b_2 \cdots b_l = e\}/K_0 
\] 
and 
\[ 
({\bf a}^{-1}(0) \cap O_\Lambda)/K_0 =  
\{(\xi_1, \xi_2, \cdots, \xi_l) \in O_\Lambda: \xi_1 + \xi_2 + \cdots + \xi_l = 0\}/K_0 
\] 
are homeomorphic. 
In particular, one is non-empty if and only if the other is. 
\end{cor} 
 
\noindent 
{\bf Proof.} Let $L: (A_0N_0)^l \to (\fp_0)^l$ be the diffeomorphism in Theorem  
\ref{thm_main-intro}. 
Then $L$ induces a homeomorphism $L: {\bf m}^{-1}(e) \to 
{\bf a}^{-1}(0)$. Since $L(D_\Lambda) = O_\Lambda$ and $L$ is $K_0$--equivariant,  
it induces a homeomorphism   
from $({\bf m}^{-1}(e) \cap D_\Lambda)/K_0$ to $({\bf a}^{-1}(0) \cap O_\Lambda)/K_0$. 
\qed 
 
\begin{rem} 
\label{rem_twisted} 
{\em 
Equip $X_0 = G_0/K_0$ with the $G_0$--invariant Riemannian structure defined 
by the restriction of the Killing form of $\fg_0$ on $\fp_0$. For $x_1, x_2 \in X_0$, let  
$\overline{x_1x_2}$ be the unique geodesic in $X_0$ connecting $x_1$ and $x_2$.  
Then there is a unique $\lambda \in \fa_{0}^{+}$ such that 
$g \cdot x_1 = \ast$ and $g \cdot x_2 = (\exp \lambda) \cdot \ast$ for some 
$g \in G_0$, where $\ast = e K_0 \in X_0$ is the base point. The element $\lambda \in \fa_{0}^{+}$ 
is called \cite{kumar-millson} the $\fa_{0}^{+}$--length of $\overline{x_1x_2}$.  
Representing the vertices of an $\ast$--based $l$--gon in $X_0$ by 
\[ 
(\ast, \, \, b_1 \cdot \ast, \, \, b_1b_2 \cdot \ast, \, \, \cdots, \, \,  
b_1b_2 \cdots b_l \cdot \ast) 
\] 
for some $b_1, b_2, \cdots, b_l \in A_0N_0$, we can regard $(b_1, b_2, \cdots, b_l)$ as the set of edges of  
the $l$--gon. Then for $\Lambda \in (\fa_{0}^{+})^l$, the space 
\[ 
\{(b_1, b_2, \cdots, b_l) \in D_\Lambda: \, b_1b_2 \cdots b_l = e\}/K_0 
\] 
can be identified with the space of $G_0$--equivalence classes 
of {\it closed} $l$--gons in $X_0$ with fixed  
``side length" $\Lambda$. Similarly, the space 
\[ 
\{(\xi_1, \xi_2, \cdots, \xi_l) \in O_\Lambda: \xi_1 + \xi_2 + \cdots + \xi_l = 0\}/K_0 
\] 
can be identified with the space of equivalent $l$--gons with fixed side length 
in the ``infinitesimal symmetric 
space" $\fp_0$. See \cite{kumar-millson} for details. 
%We give another interpretation of the twisted diagonal action on $K_0$ on 
%$(A_0N_0)^l$. Consider the action of $(K_{0})^{l}$ on $(G_{0})^{l}$ by 
%\[ 
%(k_1,k_2, \cdots, k_l)\cdot 
%(g_1, g_2, \cdots , g_n) := (k_1g_1{k_2}^{-1}, \, k_2g_2{k_3}^{-1}, \,  
%\cdots,\, k_{l-1}g_{l-1}k_{l}^{-1}, \, k_lg_l), 
%\] 
%where  $k_i \in K_0$ and $g_i \in G_0$ for $1 \leq i \leq l$. Denote by 
%\[ 
%Z = (K_{0})^{l}\backslash (G_{0})^{l}=\{ 
%[g_1, g_2, \cdots, g_l]: (g_1, g_2 \cdots, g_l) \in (G_0)^l\}  
%\] 
%the quotient space. Then $K_0$ acts on $Z$ by  
%\[ 
%k \cdot [g_1, g_2, \cdots, g_l] = [g_1, g_2, \cdots, g_l k^{-1}], 
%\] 
%and the map 
%\[ 
%{\bf m}_Z: \, \, Z \lrw K_0\backslash G_0: \, \, [g_1, g_2, \cdots, g_l] \Map K_0  
%(g_1g_2 \cdots g_l) 
%\] 
%is $K_0$-equivariant, where $K_0$ acts on $K_0 \backslash G_0$ by (\ref{eq_left-K0-action}). 
% It is easy to see that 
%the inclusion $(A_0N_0)^l\hookrightarrow (G_{0})^{l}$ induces a diffeomorphism  
%from $(A_0N_0)^l$ to $Z$. Under this diffeomorphism, the $K_0$--action on $Z$ becomes 
%the $K_0$--action on $(A_0N_0)^l$ given in (\ref{eq_twisted-action-K0}), and the  
%$K_0$--equivariant map ${\bf m}_Z: Z \to K_0 \backslash G_0 \cong A_0N_0$ becomes  
%the map ${\bf m}: (A_0N_0)^l \to A_0N_0$.  
Using the right $A_0N_0$--action on $K_0$ 
given by 
\begin{equation} 
\label{eq_AN-on-K-1} 
k^b :=  q( kb), \hs \hs b \in A_0N_0, \, k\in K_0, 
\end{equation} 
where $q: G_0 \to K_0$ is the projection $q(b_1k_1) = k_1$ for $b_1  
\in A_0N_0$ and $k_1 \in K_0$, 
it is easy to see that ${{\mathcal T}}_k$ is also given by 
\begin{equation} 
\label{eq_twisted-action-K0} 
{{\mathcal T}}_k (b_1, b_2, \cdots, b_l) :=  
(k_1 \cdot b_1, \, \, k_2 \cdot b_2, \, \cdots,\,  
k_l \cdot b_l), 
\end{equation} 
where $k_1 = k, k_{j} =k^{b_1b_2 \cdots b_{j-1}}$  
for $2 \leq j \leq l$. In the Appendix, we will see that this formula naturally 
arises in the context of Poisson Lie group actions. 
} 
\end{rem}

\begin{rem} 
\label{rem_singular} 
{\em 
%For $\Lambda = (\lambda_1, \lambda_2, \cdots, \lambda_l) \in (\fa_0)^l$, let 
%$a_j =\exp \lambda_j \in A_0$ and let $C_{a_j} = K_0a_jK_0 \subset G_0$ for  
%$1 \leq j \leq l$. Clearly $D_{a_j} \subset C_{a_j}$ for  
%each $j$. Let 
%$C_\Lambda = C_{a_1} \times C_{a_2} \times \cdots \times C_{a_l}$. Consider the  
%map  
%\[ 
%{\bf M}_{0}: \, \, (G_0)^l \lrw G_0: \, \,  
%(g_1, g_2, \cdots, g_l) \Map g_1g_2 \cdots g_l,  
%\] 
%and consider the action of $(K_0)^l$ on $(G_0)^l$ by 
%\[ 
%(k_1, k_2, \cdots, k_l) \cdot (g_1, g_2, \cdots, g_l) :=  
%(k_1g_1{k_2}^{-1}, \, k_2g_2{k_3}^{-1}, \,  
%\cdots,\, k_{l-1}g_{l-1}k_{l}^{-1}, \, k_lg_lk_{1}^{-1}). 
%\] 
%Then it is easy to show (Lemma 4.1 of \cite{amw:linear}) that 
%the natural map  
%\[ 
%({\bf m}_{0}^{-1}(e) \cap D_\Lambda)/K_0 \lrw ({\bf M}_{0}^{-1}(e) \cap C_\Lambda)/(K_0)^l 
%\] 
%is a homeomorphism. Thus by Thompson's conjecture, we have a homeomorphism 
%\[ 
%({\bf a}_{0}^{-1}(0) \cap O_\Lambda)/K_0 \, \cong \, ({\bf M}_{0}^{-1}(e) \cap C_\Lambda)/(K_0)^l. 
%\] 
When $G_0 = GL(n,\R)$, recall that the singular values of $g \in G_0$ are by definition 
the eigenvalues of $\sqrt{gg^t}$.  Thompson's conjecture for $GL(n, \R)$ says that, for any 
collection $(\lambda_1, \lambda_2, \cdots, \lambda_l)$ of real diagonal matrices, the 
following two statements are equivalent (see Section 4.2 of \cite{amw:linear}): 
 
1) there exist   
matrices $g_j \in GL(n,\R)$ whose singular values are entries of $a_j = \exp(\lambda_j)$ 
and $g_1g_2 \cdots g_l = e$; 
 
2) there exist  symmetric matrices $\xi_j$ whose eigenvalues are entries of $\lambda_j$ 
and such that $\xi_1 + \xi_2 + \cdots + \xi_l = 0$. 
} 
\end{rem} 
 
\begin{rem} 
{\em 
By a theorem of  
O'Shea and Sjamaar \cite{o-reyer:pair}, the set ${\bf a}_{0}^{-1}(0)\cap O_\Lambda$ is 
non-empty if and only if $\Lambda =(\lambda_1, \lambda_2, \cdots, \lambda_l)$ lies in a  
certain polyhedral cone in $(\fa_{0}^{+})^l$. There have been intensive research activities 
on  
the inequalities that describe this polyhedral cone. We refer to  
\cite{kumar-millson} for explicit examples of these polyhedral cones when 
$X_0$ has rank $3$ and to  
\cite{fulton:survey} \cite{KML:1} and \cite{KML:2} for an account on  
the history and the connections 
between this problem and others fields such as Schubert calculus, representation theory,  
symmetric spaces, and integrable systems. See also Remark  
\ref{rem_reyer}. 
} 
\end{rem}

\begin{rem} 
\label{rem_adjoint-type-enough} 
{\rm Finally, we remark that it is enough to prove Theorem \ref{thm_main-intro} for  
$G$ with trivial center. Indeed, let $Z$ be the center of $G$, 
let $Z_0=G_0\cap Z$, and let 
 $G_{0}^{\prime} = G_0/Z_0$ 
 with Lie algebra $\fg_{0}^{\prime} = \fg_0 /\fz_0$, where the Lie algebra
 $\fz_0$ of $Z_0$ is the center of $\fg_0$. Let $j: G_0 \to G_{0}^{\prime}$ 
and $\fg_0 \to \fg_{0}^{\prime}$ be the natural projections. Let 
$\fk_{0}^{\prime} = j(\fk_0)$ and $\fp_{0}^{\prime} = j(\fp_0)$,  
and let $K_{0}^{\prime} = j(K_0)$,  
$A_{0}^{\prime} = j(A_0)$,  and 
$N_{0}^{\prime} = j(N_0)$. By Corollary 1.3 from \cite{knapp:examples}, $Z_0 = 
(K_0 \cap Z_0) (A_0 \cap Z_0)$, and $\exp:\fz_0 \cap \fa_0 
\to Z_0 \cap A_0$ is an isomorphism. 
Thus $G_{0}^{\prime} = K_{0}^{\prime} A_{0}^{\prime} 
N_{0}^{\prime}$ is an Iwasawa decomposition for $G_{0}^{\prime}$, and  
$\fg_{0}^{\prime} = \fk_{0}^{\prime} + \fp_{0}^{\prime}$ is a Cartan decomposition  
for $\fg_{0}^{\prime}$. Moreover, $\fp_{0} \cong \fp_{0}^{\prime} 
\oplus (\fa_0 \cap \fz_0)$ and $A_0 \cong  
A_{0}^{\prime} \times (A_0 \cap Z_0)$ and $N_{0}^{\prime} \cong N_0$.  
Let  
\beqa 
{\bf a}: & & (\fp_{0}^{\prime})^l \lrw \fp_{0}^{\prime},  \hs \hs 
{\bf m}: \, \, (A_{0}^{\prime}N_{0}^{\prime})^l \lrw A_0N_0, \\ 
{\bf a}: & & (\fa_0 \cap \fz_0)^l \lrw \fa_0 \cap \fz_0,  \hs \hs 
{\bf m}: \, \, (A_{0} \cap Z_0)^l \lrw A_0 \cap Z_0 
\eeqa 
be respectively the addition and multiplication maps. If $L^\prime: (A_{0}^{\prime} 
N_{0}^{\prime})^l  \to (\fp_{0}^{\prime})^l$ is a diffeomorphism satisfying 
the requirements in Theorem \ref{thm_main-intro} for the group $G_{0}^{\prime}$, then 
$L = (L^\prime, (\log)^l)$ 
will be a diffeomorphism from $(A_0N_0)^l$ to $(\fp_0)^l$ satisfying  
the requirements in Theorem \ref{thm_main-intro} for the group $G_{0}$ , where we use the  
obvious identifications between  
$(A_0N_0)^l \cong (A_{0}^{\prime} 
N_{0}^{\prime})^l \times (A_0 \cap Z_0)^l$ and $(\fp_0)^l \cong (\fp_{0}^{\prime})^l 
\times (\fa_0 \cap \fz_0)^l$. 
} 
\end{rem}

\section{Inner classes of real forms and quasi-split real forms} 
\label{sec_inner} 
 
Let $\fg$ be a semi-simple complex Lie algebra.  
Recall that real forms of $\fg$ are in one to one correspondence  
with complex conjugate 
linear involutive automorphisms of $\fg$. For such an involution $\tau$,  
the corresponding 
real form is the fixed point set $\fg^\tau$ of $\tau$. We will refer to 
both $\fg^\tau$ and $\tau$ as the real form. Throughout this paper, if $V$ is a  
set and $\sigma$ in an involution on $V$, we will use $V^\sigma$ to denote 
the set of $\sigma$-fixed points in $V$. Let 
$G$ be the adjoint group of $\fg$. 
 
% Our reference for this section is 
%\cite{a-v:l}. 

%\begin{dfn} (Definitions 6.10 and 6.13 of \cite{a-v:l}) 
\begin{dfn} (Definitions 2.4 and 2.6 of \cite{abv:l}) 
\label{dfn_steinberg} 
{\em Two real forms $\tau_1$ and $\tau_2$ of $\fg$ 
are said to be 
inner to each other if there exists $g \in G$ such that 
$\tau_1 = \Ad_g \tau_2$. 
A real form $\tau$ of $\fg$ is said to be {\it quasi-split} if there exists 
a Borel subalgebra $\fb$ of $\fg$ such that $\tau(\fb) = \fb$. 
} 
\end{dfn} 
 
Inner classes 
of real forms are classified by involutive automorphisms of the 
Dynkin diagram $D(\fg)$ of $\fg$. Indeed,  
let $\Autg$ be the group of complex linear automorphisms of $\fg.$ 
Its identity component is the adjoint group $G$. 
Let $\Outs$ be the automorphism 
group of the Dynkin diagram $D(\fg)$ of $\fg$.  
%It is well-known (see, for example, Section 2.15 of  
%\cite{samelson:notes}) 
%that 
There 
is a split short exact sequence (\cite{abv:l}, Proposition 2.11), 
\begin{equation} 
\label{eq_short-exact} 
1 \, \longrightarrow G \, \longrightarrow \, \As \, 
\stackrel{\varrho}{\longrightarrow} \, \Outs \, \longrightarrow \, 1. 
\end{equation} 
Denote by ${\bf R}$ the set of all real forms of $\fg$. Let $\theta$ be 
any compact real form of $\fg$. 
Define 
\begin{equation} 
\label{eq_psi} 
\varpi: \, {\bf R} \, \longrightarrow \, \Outs: \, 
\varpi(\tau) \, = \, \varrho(\tau \theta). 
\end{equation}  
Then $\varpi(\tau)^2 = 1$ 
for each $\tau$, and $\tau_1$ and $\tau_2$ 
are inner to each other if and only if $\varpi(\tau_1) = \varpi(\tau_2)$. 
Conversely, for every involutive $d \in {\rm Aut}_{D({\mathfrak g})}$, we can 
construct $\gamma_d \in {\rm Aut}_{{\mathfrak g}}$ such that  
$\tau:=\gamma_d \theta$ is a real form with $\varpi(\tau) = d$  
(see (\ref{eq_gamma-d}) below). 
Thus the map $\varpi$  
gives a bijection between inner classes of real forms of $\fg$ and  
involutive elements in $\Outs$ (Proposition 2.12 of  
\cite{abv:l}).   
 
\begin{dfn} 
\label{dfn_class-d} 
{\em 
Let $d$ be an involutive automorphism of the Dynkin 
diagram $D(\fg)$ of $\fg$. We say that a real form $\tau$ of $\fg$ 
is {\it of inner class $d$} or {\it in the  $d$--inner class} if $ 
\varpi(\tau) = d$. 
} 
\end{dfn} 
 
By Proposition 2.7 of \cite{abv:l}, every inner class of real forms of $\fg$ 
contains a quasi-split real form that is 
unique up to $G$--conjugacy. 
In the following, for each involutive $d \in \Outs$, we will construct  
an explicit quasi-split real form 
$\tau_d$ in the $d$--inner class. We will then show that, up 
to $G$--conjugacy, every  real form in the $d$--inner class  is 
of the form 
$\tau = \Ad_{\dot{w}_0} \tau_d$, where $w_0$ is a certain Weyl group element 
and $\dot{w}_0$ a representative of $w_0$ in $G$. 
  We 
 first fix once and for all the following data for $\fg$:  
 
Let $\fh$ be a Cartan subalgebra of $\fg$,  and  
let $\Delta$ be the corresponding root system. Fix a choice of positive  
roots $\Delta^+$in $\Delta$,  and let $\Sigma$ be the basis of simple roots.  
Let $\ll\cdot,\cdot \gg$ be the Killing form on $\fg$. 
For each $\alpha \in \Delta$, 
let 
$E_{\pm \alpha}$ be root vectors such that $[E_\alpha, E_{-\alpha}] = 
H_\alpha$ for all $\alpha \in \Delta^+$, 
where $H_\alpha$ is the unique element of $\fh$ defined by $\ll H, 
H_\alpha\gg =\alpha(H)$ for all $H\in\fh$, and the numbers  
$m_{\alpha,\beta}$ 
for $\alpha, \beta \in \Delta$ defined for $\alpha\ne -\beta$ by 
\[ 
\begin{array}{cccc} 
[E_\alpha, E_\beta] & = & m_{\alpha, \beta}E_{\alpha+\beta} & \mbox{ if } 
\alpha+\beta\in \Delta \\ 
{} & = & 0 & \mbox{otherwise}     
\end{array} 
\] 
 satisfy $m_{-\alpha, 
-\beta}=-m_{\alpha, \beta}$. We will refer to the set  
$\{E_\alpha, E_{-\alpha}: \alpha\in \Delta^+\}$ as (part of) a  
Weyl basis. 
Using a Weyl basis $\{E_\alpha, E_{-\alpha}: \alpha\in \Delta^+\}$, we can  
define a compact real form $\fk$ of $\fg$ as 
\begin{equation} 
\label{eq_fk} 
\fk =  
{\rm span}_{{\mathbb R}} \{ \i H_\alpha, X_\alpha:=E_\alpha-E_{-\alpha},  
Y_\alpha:=\i(E_{\alpha}+E_{-\alpha}): \, \alpha \in  
\Delta^+\}. 
\end{equation} 
Let $\theta$ be the complex conjugation of $\fg$ 
defining $\fk$.  We also  define a split real form 
$\eta_0$ of $\fg$ by setting 
$\eta_0 |_{{\mathfrak a}} = id$, and $\eta_0(\ea) = \ea$ for every $\alpha \in \Delta$. 
Clearly $\theta \eta_0 = \eta_0 \theta$. The inner class for 
 $\eta_0$ is easily  
seen to be the automorphism of the simple roots given by 
$-w^0$, where $w^0$ is the longest element in the Weyl group $W$ 
of $(\fg, \fh)$.

An explicit splitting of the 
short exact sequence (\ref{eq_short-exact}) can be constructed using the 
Weyl basis. Indeed,  
for any $d \in \Outs$, define $\gamma_d \in \Autg$ by requiring 
\begin{equation} 
\label{eq_gamma-d} 
\gamma_d(H_\alpha) = H_{d\alpha}, \hs {\rm and} \hs  
\gamma_d (E_\alpha ) = E_{d \alpha} 
\end{equation} 
for each simple root $\alpha$. Then $d \mapsto \gamma_d$ is a group 
homomorphism from $\Outs$ to $\Autg$ and is a section of $\varrho$ in  
(\ref{eq_short-exact}).  
Moreover,  
every $\gamma_d$ commutes with both $\theta$ and $\eta_0$ because they commute  
on a set of 
generators of $\fg$.  
 
\begin{lem} 
\label{lem_dfn-tau-d} 
For an involutive element  $d \in \Outs$,   
let $\gamma_{-w^0d}\in {\rm Aut}_{{\mathfrak g}}$ 
be the lifting of $-w^0d \in {\rm Aut}_{D({\mathfrak g})}$ 
as defined in (\ref{eq_gamma-d}). Define 
\begin{equation} 
\label{eq_tau-d}  
\tau_d = \eta_0 \gamma_{-w^0d}. 
\end{equation} 
Then $\tau_d$ is a quasi-split real form of $\fg$ in the $d$--inner class. 
\end{lem} 
 
\noindent 
{\bf Proof.} We know that $(\tau_d)^2 = 1$ because $-w^0 \in \Outs$ is in the center. 
Since $\tau_d$ maps every positive root vector to another positive root vector,  
it is a quasi-split real form. Finally, since $\varpi(\tau_d \theta)  
=\varrho(\gamma_{-w^0} \gamma_{-w^0d}) = d$, we see that $\tau_d$ is  
in the $d$--inner class. 
\qed 
  
To relate an arbitrary real form  
in the $d$--inner class with the quasi-split real form $\tau_d$,  
we recall some definitions from  
\cite{araki:satake}. Note first that $\fg= 
\fk + \fa + \fn$ is an Iwasawa decomposition, where 
$\fa = {\rm span}_{{\mathbb R}}\{H_\alpha: \alpha \in \Delta\}$ 
and $\fn = \sum_{\alpha \in \Delta^+} \fg_\alpha$.

\begin{dfn} 
\label{dfn_iwasawa} 
{\em 
A real form $\tau$ of $\fg$ is said to be {\it normally related to $(\fk, \fa)$} 
and {\it compatible with $\Delta^+$} if  
 
1) $\tau \theta = \theta \tau$, and $\tau(\fh) = \fh$; 
 
2) $\fa^\tau \subset (\sqrt{-1}\fk)^\tau$ is maximal abelian in $(\sqrt{-1}\fk)^\tau$; 
 
3) if $\alpha \in \Delta^+$ is such that $\alpha|_{{\mathfrak a}^{\tau}}  
\neq 0$, then $\tau(\alpha) \in \Delta^+$, 
where $\tau(\alpha) \in \fa^*$ is defined by $\tau(\alpha)(\lambda) =  
\alpha(\tau(\lambda))$ for $\lambda \in \fa$.

\noindent 
We will call a real form with Properties 1)-3) an {\it Iwasawa real form} 
relative to $(\fk, \fa, \fn)$. 
} 
\end{dfn} 
 
\begin{rem} 
\label{rem_iwasawa} 
{\em 
Once 1) and 2) in Definition \ref{dfn_iwasawa} are satisfied, 3) is equivalent to 
the set $r(\Delta^+) \backslash \{0\} \subset \Delta_{{\rm res}}$ being a choice 
of positive roots for the restricted root system $\Delta_{{\rm res}}$ of 
$(\fg^\tau, \fa^\tau)$, where $r$ is the map dual to the inclusion 
$\fa^\tau \hookrightarrow \fa$.  
%It is also easy to see that 
If $\tau$ is an Iwasawa real form relative to $(\fk, \fa, \fn)$, then so 
is $\Ad_t \tau \Ad_{t}^{-1}$ for any $t \in \exp (i\fa)$. 
} 
\end{rem}

\begin{prop} 
\label{prop_tau-s} 1) Every real form of $\fg$ is  
conjugate by an element in $G$ to a real form that is Iwasawa  
relative to $(\fk, \fa, \fn)$; 
 
2) Suppose that $\tau$ is an Iwasawa real form 
relative to $(\fk, \fa, \fn)$ and suppose that $\tau$ is in the $d$--inner class. 
Let $w_0$ be the longest element of the subgroup of $W$ generated by  
the reflections corresponding to roots in the set 
 \[ 
\Delta_{0} = \{\alpha \in \Delta: \alpha|_{{\mathfrak a}^\tau} = 0\}=\{\alpha \in \Delta: 
\, \tau(\alpha) = -\alpha\}. 
\] 
Then 
there is a representative $\dot{w}_0$ of $w_0$ in $G$ such that  
\begin{equation} 
\label{eq_tau} 
\tau = \Ad_{\dot{w}_0} \tau_d. 
\end{equation} 
\end{prop} 
  
\noindent 
{\bf Proof.} Statement 1) follows from Proposition 1.2, Section 2.8, and Corollary 2.5  
 of \cite{araki:satake}.  
 
 2) Assume now that $\tau$ is an Iwasawa real form 
relative to $(\fk, \fa, \fn)$ and that $\tau$ is in the $d$--inner class.  
Consider $\tau_d \tau$. Since both $\tau$ and $\tau_d$ are  
in the $d$-inner class, $\tau_d \tau = \Ad_g$ for some $g \in G$. 
Since both $\tau$ and $\tau_d$ leave $\fh$ 
invariant, the element $g$ represents an element in the Weyl group $W$. 
Let $\Sigma_0 = \Sigma \cap \Delta_{0}$. By Section 2.8 of 
\cite{araki:satake}, $\alpha \in \Delta_0$ if and only if  
$\alpha$ is in the linear span of $\Sigma_0$.  
For every $\alpha \in \Sigma_0$, we have $\tau_d \tau(\alpha) = -\tau_d ( 
\alpha) \in -\Delta^+$, and for every $\alpha \in \Delta^+ - \Delta_0$,  
since $\tau(\alpha) \in \Delta^+$, we have $(\tau_d \tau)(\alpha) \in 
\Delta^+$. Thus $g$ represents the element $w_0= (w_0)^{-1}$. 
\qed

\begin{rem} 
\label{rem_satake-diagram} 
{\em 
Recall from  
\cite{araki:satake} that the Satake  
diagram of $\tau$ is the Dynkin diagram of $\fg$ with simple roots in  
$\Sigma_0$ painted black, simple roots in $\Sigma - \Sigma_0$ painted  
white, and a two sided arrow drawn between two white  
simple roots $\alpha$ and $\alpha^\prime$ if $\tau(\alpha) =  
\alpha^\prime + \beta$ for some $\beta \in \Delta_0$. From  
(\ref{eq_tau}) we see that $\alpha^\prime = -w^0d (\alpha)$ if 
$\alpha$ is a white simple root. Conversely, given a Satake 
diagram for a real form of $\fg$, let 
$c \in {\rm Aut}_{D({\mathfrak g})}$ be defined by  
\begin{equation} 
\label{eq_cs} 
c(\alpha) = \begin{cases} -w_0 \alpha & {\rm if}  \hspace{.05in} 
 \alpha \hspace{.05in}  
{\rm is}  \hspace{.05in} {\rm black} \\ 
\alpha^\prime & {\rm if}  \hspace{.05in} \alpha  \hspace{.05in} 
{\rm is}  \hspace{.05in} {\rm white}, \end{cases} 
\end{equation} 
where $w_0$ is the longest element 
in the subgroup of the Weyl group 
of $(\fg, \fh)$ generated by the black dots in the Satake diagram, and  
$\alpha \mapsto \alpha^\prime$ is the order $2$ involution on the set of  
white dots in the diagram.  
Then $c$ is involutive, and the inner class of the real form is $d =-w^0 c$. 
} 
\end{rem} 
   
We now return to the real form  $\tau$ in Proposition \ref{prop_tau-s}.  
Set 
\[ 
\Delta_{1}^{+} = \{\alpha \in \Delta^+: \alpha|_{{\mathfrak a}^\tau} \neq 0\}, 
\hs 
(\Delta_{1}^{+})^\prime = \Delta^+ \cap \Delta_0 = \{\alpha \in \Delta^+:  
\alpha|_{{\mathfrak a}^\tau} = 0\}. 
\] 
Then $\tau(\Delta_{1}^{+}) \subset \Delta_{1}^{+}$. 
Set 
\begin{equation} 
\label{eq_fn-1-2} 
\fn_{1} = \sum_{\alpha \in \Delta_{1}^{+}} \fg_\alpha, \hs 
\fn_{1}^{\prime} = \sum_{\alpha \in (\Delta_{1}^{+})^\prime} \fg_\alpha. 
\end{equation} 
Then $\fn_1$ is $\tau$-invariant, and, since there are no non-compact  
imaginary roots for the Cartan subalgebra $\fh^\tau$ of $\fg^\tau$, 
we have $\tau|_{{\mathfrak n}_{1}^{\prime}} = \theta|_{{\mathfrak n}_{1}^{\prime}}$  
(Proposition 1.1 of 
\cite{araki:satake}). 
Since the restriction of $\Delta^+$ to $\fa^\tau$ gives a choice of  
positive restricted roots for $(\fg^\tau, \fa^\tau)$, we know that 
\begin{equation} 
\label{eq_iwasawa-tau-s} 
\fg^{\tau} = \fk^{\tau} + \fa^{\tau} + (\fn_1)^{\tau} 
\end{equation} 
is an Iwasawa decomposition of $\fg^{\tau}$. 
 
\section{The proof of Theorem \ref{thm_main-intro}} 
\label{sec_the-proof} 
 
By Remark \ref{rem_adjoint-type-enough}, it is enough to 
prove Theorem \ref{thm_main-intro} when $G$ of adjoint type.  
Let $\fg_0$ be a real form of $\fg$, the Lie algebra of $G$.
 Then by Proposition \ref{prop_tau-s}, we can assume 
that $\fg_0=\fg^\tau$, where $\tau$ is the involution on $\fg$ given by 
(\ref{eq_tau}), and $d$ is the inner class of $\fg_0$. 
 The lifting of 
$\tau$ to $G$ will also be denoted by $\tau$.  
Let $G_0$ contain
 the connected component of the identity of the subgroup $G^{\tau}$.  
%Since an element $g$ in the center of $G_0$ acts   
%trivially on $\fg = \fg_0 + i\fg_0$ for the adjoint action,
% $g \in G $ must be the  
%identity element. Thus $G_0$ has 
%trivial center.
 In this section, we will prove Theorem \ref{thm_main-intro} 
for $G_0$. 
 
The first step in out proof of Theorem \ref{thm_main-intro} for $G_0$ 
is to realize the various objects  
associated to $G_0$ as fixed point sets of involutions on  
the corresponding objects related to $G$. 
We will then apply a theorem of Alekseev-Meinrenken-Woodward, stated as 
Theorem \ref{thm_amw} below, whose proof using Poisson geometry will be outlined 
in Section \ref{sec_amw} the Appendix. 
 
We will keep all the notation from Section \ref{sec_inner}.  
In  
particular, set 
\[ 
\fk_0 = \fk^{\tau},  \hs \fp_0 = (\sqrt{-1}\fk)^\tau, \hs  
\fa_0 = \fa^{\tau}, \hs {\rm and} \hs \fn_0 = (\fn_1)^{\tau}. 
\] 
Then $\fg_0 = \fk_0 + \fp_0$ is a Cartan decomposition of  $\fg_0$, 
and $\fg_0 = \fk_0 + \fa_0 + \fn_0$ an Iwasawa decomposition of $\fg_0$.  
Let $K$ be the connected subgroup of $G$ with Lie algebra $\fk$. Let 
$K_0 = K \cap G_0$, $K^{\tau}=K\cap G^{\tau}$, and let 
\[ 
N_1 = \exp(\fn_1), \hs A_0 =  \exp(\fa_0), \hs {\rm and} \hs  
N_0 = N_1 \cap G^{\tau} = \exp(\fn_0). 
\]

\begin{lem} 
\label{lem_iwasawa-tau-s} 
$G_0 = K_0 A_0 N_0$ is an Iwasawa decomposition of $G_0$, and  
$G^\tau = K^\tau A_0N_0$ is  
an Iwasawa decomposition of $G^\tau$. 
\end{lem}

\noindent 
{\bf Proof.} The statements follow from  Lemma 2.1.7 in \cite{wallach:book} 
and the facts that $K_0$ and $K^\tau$ are  maximally 
compact subgroups of $G_0$ and of $G^\tau$ respectively. 
\qed 
  
Let $A = \exp \fa$ and $N = \exp \fn$ so that $G=KAN$ 
is an Iwasawa decomposition of $G$. 
We will now identify 
$A_0N_0$ with the fixed point set of an involution on $AN$, although we note that  
the involution $\tau$ does not leave $AN$ invariant unless $\tau = \tau_d$. 
Define 
\begin{equation} 
\label{eq_sigma-s} 
\sigma: \, AN \lrw AN: \, \sigma (b) = p (\tau(b)), 
\end{equation} 
where $p: G =KAN \to AN$ is the projection $b_1k_1 \mapsto b_1$ for $k_1 \in K$ and 
$b_1 \in AN$. 
Recall from (\ref{eq_tau}) that $\tau = {\rm Ad}_{\dot{w}_0} \tau_d$ on $\fg$. 
Since both $\tau$ and $\tau_d$ commute with $\theta$, the element  
$\dot{w}_0$ is in $K$. As for the case of $G_0$, we can use the Iwasawa 
decomposition $G = KAN$ to define a left action of $K$ on $AN \cong  
G/K$ by 
\begin{equation} 
\label{eq_K-on-AN} 
k \cdot b := p(kb), \hs \hs k \in K, \, \, b \in AN. 
\end{equation} 
Then $\sigma: AN \to AN$ is also given by 
\begin{equation} 
\label{eq_sigma-tau-d} 
\sigma (b) = \dot{w}_0 \cdot \tau_d(b), \hs {\rm for} \hs b \in AN. 
\end{equation}

\begin{lem} 
\label{lem_sigma-s-involution} 
$\sigma: AN \to AN$ is an involution, and  
$A_0N_0 = (AN)^{\sigma}$, the fixed 
point set of $\sigma$ in $AN$. 
\end{lem} 
 
\noindent 
{\bf Proof.} The fact that $\sigma^{2} = 1$ follows from 
 the fact that $\dot{w}_0 \tau_d(\dot{w}_0) = 1$. 
Recall that $A_0N_0$ is the fixed point set of $\tau$ in $AN_1$, where 
$N_1 = \exp(\fn_1)$ and $\fn_1 \subset \fn$ is given in (\ref{eq_fn-1-2}). 
Since $\sigma$ coincides with $\tau$ on $AN_1$, we have  
$A_0N_0 \subset (AN)^{\sigma}$. Suppose now that  
$b \in AN$ is such that $\sigma(b) = b$. Then there exists  
$k \in K$ such that $\tau(b) = bk$. By Proposition 7.1.3 in 
\cite{schli:symmetric}, we can decompose $b$ as $b = gak_1$ for some $k_1 \in K,  
a \in A$ and $g \in G^{\tau}$. Then $\tau(b) = g \tau(a)\tau(k_1)$, 
and thus $g \tau(a)\tau(k_1) = gak_1k$, from which it 
follows that $\tau(a) = a$ and 
$\tau(k_1) = k_1k$. Thus $k =k_{1}^{-1}\tau(k_1)$, and hence $\tau(bk_{1}^{-1}) =b k_{1}^{-1}$.  
Write  $bk_{1}^{-1}=b_2k_2$ with 
$k_2\in K^{\tau}$ and $b_2\in A_0N_0$ using the Iwasawa decomposition 
of $G^\tau$. It follows then that $k_2=k_{1}^{-1}$ and $b_2=b$, 
so $b\in A_0N_0$. 
\qed 
 
Let now $l \geq 1$ be an integer. As in Section \ref{sec_the-conjecture}, we have the 
twisted diagonal action $k \mapsto {{\mathcal T}}_k$ of $K$ on $(AN)^l$ given by 
\begin{equation} 
\label{eq_twisted-action-K} 
{{\mathcal T}}_k = \nu^{-1} \circ \delta_k \circ \nu: \, \, (AN)^l \lrw (AN)^l, 
\end{equation} 
where $\nu: (AN)^l \to (AN)^l$ is as in (\ref{eq_nu}) with $A_0N_0$ replaced by $AN$, and  
$\delta_k$ denotes the diagonal action of $k \in K$ on $AN$. Set 
$(\tau_d)^l = (\tau_d, \tau_d, \cdots,  
\tau_d): \,  (AN)^l \to (AN)^l.$

\begin{lem} 
\label{lem_product} 
For an integer $l \geq 1$, define 
\[ 
\sigma^{(l)} = {{\mathcal T}}_{\dot{w}_0} \circ (\tau_d)^{l}: \, \, (AN)^l \lrw (AN)^l. 
\] 
Then $\sigma^{(l)}$ is an involution, and the fixed point set of $\sigma^{(l)}$ is 
$(A_0N_0)^l$. 
\end{lem} 
  
\noindent 
{\bf Proof.} Let $\sigma^l = (\sigma, \sigma, \cdots, \sigma): (AN)^l 
\to (AN)^l$, where $\sigma: AN \to AN$ is as in Lemma \ref{lem_sigma-s-involution}. 
Since $\tau_d$ is a group automorphism of $AN$, we have 
\[ 
\sigma^{(l)} = \nu^{-1} \circ \delta_{\dot{w}_0} \circ \nu \circ (\tau_d)^l =  
\nu^{-1} \circ (\delta_{\dot{w}_0} \circ (\tau_d)^l ) \circ \nu =  
\nu^{-1} \circ \sigma^l \circ \nu. 
\] 
Thus  
$(\sigma^{(l)})^2 = 1$. Moreover,  let $b=(b_1, b_2, \cdots, b_l) \in (AN)^l$, and let 
$b^\prime= \nu(b)$. Then $\sigma^{(l)}(b) = b$ 
if and only if $\sigma^l(b^\prime) = b^\prime$, which in turn is equivalent to $b_j \in A_0N_0$ for each  
$1 \leq j \leq l$ because of  Lemma 
\ref{lem_sigma-s-involution} and because of the fact that $A_0N_0$ is a subgroup of $AN$. 
\qed 
 
Let $\fp = \sqrt{-1} \fk$, so $\fg = \fk + \fp$ is a Cartan decomposition of $\fg$. 
Let $E: \fp \to AN$ be the composition of the identifications 
\begin{equation} 
\label{eq_E} 
E: \, \, \fp \stackrel{\exp}{\cong} \exp(\fp) \cong G/K \cong AN. 
\end{equation} 
Then $E$ is $K$--equivariant with respect to the action of $K$ on $\fp$ by  
conjugation and the action of $K$ on $AN$ given in  
(\ref{eq_K-on-AN}). 
 
\begin{lem} 
\label{lem_E-sigma-s} 
$E \circ (\tau|_{{\mathfrak p}})= \sigma \circ E,$ and $E_0 = E|_{{\mathfrak p}_0}: 
\fp_0 \to A_0N_0$. 
\end{lem} 
 
\noindent 
{\bf Proof.} Consider $E^{-1}: AN \to \fp$. For $g \in G$, define 
$g^* = \theta(g^{-1})$. Then for every $b \in AN$,  
$E^{-1}(b) = {1 \over 2} \log (bb^*)$ for all $b \in AN$, and thus 
\beqa 
{E^{-1}}(\sigma(b))& = &\frac{1}{2}\log(\dot{w}_0 \tau_d(b)\tau_d(b)^*{\dot{w}_0}^{-1})\\ 
& = & {\rm Ad}_{\dot{w}_0}\left(\frac{1}{2}\log( \tau_d(b)\tau_d(b)^*)\right) 
 = {\rm Ad}_{\dot{w}_0}\tau_d(E^{-1}(b)) = \tau (E^{-1}(b)). 
\eeqa 
Thus $E \circ (\tau|_{{\mathfrak p}})= \sigma \circ E,$ and $E(\fp_0) = A_0N_0$ 
by Lemma \ref{lem_sigma-s-involution}. It also follows that  
$E|_{{\mathfrak p}_0} = E_0$. 
\qed

\begin{nota} 
\label{nota_3} 
{\em 
For $\lambda \in \fa\subset \fp $, let ${{\mathcal O}}_\lambda$ 
be the $K$--orbit in $\fp$ through $\lambda$, and let ${{\mathcal D}}_a$ 
be the $K$--orbit in $AN$ through $a = \exp \lambda  \in A$.  
For $\Lambda =(\lambda_1, \lambda_2, \cdots, \lambda_l) 
\in (\fa)^l$, set 
\begin{equation} 
\label{eq_O} 
{{\mathcal O}}_\Lambda  = {{\mathcal O}}_{\lambda_1} \times {{\mathcal O}}_{\lambda_2} 
\times \cdots \times {{\mathcal O}}_{\lambda_l}, \hs {\rm and} \hs 
{{\mathcal D}}_\Lambda  = {{\mathcal D}}_{a_1} \times {{\mathcal D}}_{a_2} 
\times \cdots \times {{\mathcal D}}_{a_l}. 
\end{equation} 
} 
\end{nota} 
 
Recall from Section \ref{sec_the-conjecture} that, for $\lambda \in \fa_0$,  
$O_\lambda$ denotes the 
$K_0$--orbit in $\fp_0$ through $\lambda$, and $D_a$ 
denotes the $K_0$--orbit in $A_0N_0$ through $a = \exp \lambda \in A_0$.

\begin{lem} 
\label{lem_orbits} 
Let $\lambda \in \fa_{0}$ and $a = \exp \lambda \in A_0$. Then  
 
1). $\tau_d(\lambda) = \lambda$, so both  
${{\mathcal O}}_\lambda \subset \fp$ and ${{\mathcal D}}_a \subset AN$ 
are $\tau_d$-invariant; 
 
2). ${{\mathcal O}}_\lambda$ is $\tau$-invariant, and  
$({{\mathcal O}}_\lambda)^{\tau} = O_\lambda$; 
 
3). ${{\mathcal D}}_a$ is $\sigma$-invariant and  
$({{\mathcal D}}_a)^{\sigma}= D_a$. 
\end{lem} 
 
\noindent 
{\bf Proof.} 1). Let $\alpha$ be a simple root such that 
$\tau(\alpha) = -\alpha$. Then for any $\lambda \in \fa_0$, 
\[ 
\alpha(\lambda) = \alpha(\tau(\lambda)) = (\tau (\alpha))(\lambda)= 
-\alpha(\lambda) = 0. 
\] 
Thus $r_\alpha (\lambda) = \lambda$, where $r_\alpha$ is the reflection in  
$\fa$ 
defined by $\alpha$. Since $w_0$ is a product of such 
reflections (see Proposition \ref{prop_tau-s}), 
 we see that $w_0$ acts trivially on $\fa_0$. Thus 
every $\lambda \in \fa_0$ is also fixed by $\tau_d$. 
 
2). This is a standard fact. See, for example, Example 2.9 in \cite{o-reyer:pair}. We remark that the key point is to show that
 $({{\mathcal O}}_\lambda)^{\tau}$ is connected.
3) follows from 2) and Lemma \ref{lem_E-sigma-s}. 
\qed

We can now state the Alekseev-Meinrenken-Woodward 
theorem.  
Set  
\beqa 
& & {\bf a}: \,  
 \fp \times \fp \times \cdots \times \fp \lrw \fp: \, (x_1, x_2, \cdots, x_l) \Map x_1 + 
 \cdots + x_l,\\ 
& & {\bf m}: \, AN \times AN \times \cdots \times AN  \lrw AN: \, \,  
 \,  (b_1, b_2, \cdots, b_l) \Map b_1b_2 
 \cdots b_l. 
\eeqa 
As for the case of $G_0$, we will equip 
$\fp^l$ with the diagonal $K$--action by conjugation, and we will equip 
$(AN)^l$ with the twisted diagonal action  
${{\mathcal T}}$ given by (\ref{eq_twisted-action-K}). 
 
\begin{thm}[Alekseev-Meinrenken-Woodward] \cite{amw:linear} 
\label{thm_amw} 
For every quasi-split real form  
$\tau_d$ given in (\ref{eq_tau-d}) and for 
every integer $l \geq 1$, there is a $K$--equivariant diffeomorphism  
$L: (AN)^l \to \fp^l$ such that  
\begin{equation} 
\label{eq_conditions} 
{\bf m} = E \circ {\bf a} \circ L, \hs {\rm and} \hs 
 (\tau_d)^l \circ L = L \circ (\tau_d)^l. 
\end{equation} 
Moreover,  
$L({{\mathcal D}}_\Lambda) = {{\mathcal O}}_\Lambda$ for every $\Lambda \in \fa^l$. 
\end{thm} 
 
\begin{rem} 
\label{rem_appendix-amw} 
{\em 
Theorem \ref{thm_amw}, whose proof will be outlined in the Appendix,   
 is a consequence of a Moser Isotopy Lemma for Hamiltonian $K$--actions on Poisson 
manifolds, a rigidity theorem for such spaces.  More precisely, we will 
see in the Appendix that the map $E^{-1} \circ {\bf m}: (AN)^l \to \fp$ is a  
moment map for the 
twisted diagonal $K$--action ${{\mathcal T}}$ on $(AN)^l$ with respect to
 a   
Poisson structure $\pi_1$ on $(AN)^l$, and $(\tau_d)^l$ is an anti--Poisson involution 
for $\pi_1$. Moreover, the symplectic leaves of $\pi_1$ 
are precisely all the orbits ${{\mathcal O}}_\Lambda$ for $\Lambda \in \fa^l$. 
The quintuple  
\[ 
{{\mathcal Q}}_1 =((AN)^l, \, \, \pi_1, \, \, {{\mathcal T}}, \, \,  
E^{-1}\circ {\bf m}, \, \, (\tau_d)^l) 
\] 
will be referred to as a Hamiltonian Poisson $K$--space with anti--Poisson involution.  
In fact, ${{\mathcal Q}}_1$ 
 belongs to  
a smooth family  
\[ 
{{\mathcal Q}}_s = ((AN)^l, \, \, \pi_s, \, \, {{\mathcal T}}_s,\, \,   
E^{-1} \circ {\bf m}_s,\, \, (\tau_d)^l) 
\] 
as the case for $s = 1$, and when $s = 0$,  
${{\mathcal T}}_s$ is the  
diagonal $K$--action on $(AN)^l$ and $E^{-1} \circ {\bf m}_0 = {\bf a} \circ (E^{-1})^l$. 
The Moser Isotopy Lemma, Proposition  \ref{prop_moser} in Section \ref{sec_amw}, 
implies that ${{\mathcal Q}}_s$ is isomorphic to  
${{\mathcal Q}}_0$  by a diffeomorphism $\psi_s$ of $(AN)^l$ for every $s \in \R$. The map 
$L$ in Theorem \ref{thm_amw} is then taken to be $(E^{-1})^l \circ \psi_1$.  
} 
\end{rem} 
 
We will assume Theorem \ref{thm_amw} for now and prove  
Theorem \ref{thm_main-intro} for $G_0$.  
 
\bigskip 
\noindent 
{\bf Proof of Theorem \ref{thm_main-intro}.} Let $L: (AN)^l \to \fp^l$ be as in 
Theorem \ref{thm_amw}. Since $L$ is $K$--equivariant and 
intertwines $(\tau_d)^l: (AN)^l \to (AN)^l$ 
and $(\tau_d)^l: \fp^l \to \fp^l$, it 
also intertwines  
\[ 
\sigma^{(l)} = {{\mathcal T}}_{\dot{w}_0} \circ (\tau_d)^l: \, (AN)^l \lrw (AN)^l 
\hs {\rm and} \hs  
\tau^l = \delta_{\dot{w}_0} \circ (\tau_d)^l: \, \fp^l \lrw \fp^l. 
\] 
Thus by Lemma \ref{lem_product}, we know that $L((A_0N_0)^l) = (\fp_0)^l$. 
Denote $L|_{(A_0N_0)^l}: (A_0N_0)^l \to (\fp_0)^l$ also by $L$. Then clearly 
$L$ is $K_0$--equivariant, and since $E_0: \fp_0 \to A_0N_0$ coincides with 
the restriction of $E: \fp \to AN$ to $\fp_0$, we see that ${\bf m}  
= E_0 \circ {\bf a} \circ L$. Finally, let 
$\Lambda= (\lambda_1, \lambda_2, \cdots, \lambda_l) \in (\fa_0)^l$. Then 
by Lemma \ref{lem_orbits}, we see that 
\[ 
L(D_\Lambda) = L({{\mathcal D}}_\Lambda \cap (A_0N_0)^l) = {{\mathcal O}}_\Lambda 
\cap (\fp_0)^l = O_\Lambda. 
\] 
\qed 
 
\begin{rem} 
\label{rem_reyer} 
{\em 
Let $\fa^+ 
=\{x \in \fa: \alpha(x) \geq 0, \, \forall \alpha \in \Delta^+\}$ so that  
$\fa_{0}^{+} = \fa_0 \cap \fa^+$.  By Kirwan's convexity 
theorem \cite{o-reyer:pair}, there exists a polyhedral cone ${{\mathcal P}} \subset 
(\fa^+)^l$ such that $\Lambda\in {{\mathcal P}}$ 
if and only if ${\bf a}^{-1}(0) \cap {{\mathcal O}}_\Lambda$ is non-empty. 
 Set  
${{\mathcal P}}_0 = {{\mathcal P}} \cap (\fa_{0})^l \subset (\fa_{0}^{+})^l$. 
Then by a theorem of O'Shea-Sjamaar \cite{o-reyer:pair}, the set  
${\bf a}^{-1}(0) \cap O_\Lambda$  
is non-empty if and only if  
  $\Lambda \in {{\mathcal P}}_0$. If we use 
${{\mathcal P}}_d \subset (\fa^{\tau_d})^l$ to denote the polyhedral cone  
${{\mathcal P}} \cap (\fa^{\tau_d})^l$ for the quasi-split 
form $\tau_d$, it follows from $\fa_0 \subset \fa^{\tau_d}$ that  
\begin{equation} 
\label{eq_P} 
{{\mathcal P}}_0 = {{\mathcal P}}_d \cap (\fa_0)^l. 
\end{equation}  
A statement related to this fact is given in \cite{philip:loci}. } 
\end{rem} 
 
%We note that the compatibility of $\tau_{d}^{D}$ and $\tau_{d}^{O}$ with the 
%$K$-actions and the maps $\Psi$ and $\Phi$ can be expressed as 
%\begin{eqnarray} 
%\label{eq_tau-t} 
%\tau_{d}^{D}\circ t_k  & =& t_{\tau_d(k)}\circ \tau_{d}^{D}: \, {{\mathcal D}} \lrw 
%{{\mathcal D}}\\ 
%\tau_{d}^{O} \circ \Ad_k & = & \Ad_{\tau_d(k)} \circ \tau_{d}^{O}: \,  
%{{\mathcal O}} \lrw {{\mathcal O}}, 
%\end{eqnarray} 
%and 
%\begin{equation} 
%\label{eq_Psi-Phi-tau} 
%{{ \Phi}} \circ \tau_{d}^{O} = \tau_d \circ {{\Phi}} \hs {\rm and} \hs 
%{{\Psi}} \circ \tau_{d}^{D} = \tau_d \circ {{\Psi}}. 
%\end{equation} 

\section{Appendix: Alekseev-Meinrenken-Woodward theorem} 
\label{sec_amw} 
 
In this appendix, we 
give an outline of the proof of the 
 Alekseev-Meinrenken-Woodward theorem, stated as Theorem \ref{thm_amw} in 
this paper.  
Theorem 3.5 of \cite{amw:linear} shows the existence of a diffeomorphism $L_\Lambda:  
{{\mathcal D}}_\Lambda\to 
{{\mathcal O}}_\Lambda$ satisfying (\ref{eq_conditions}) 
for each $\Lambda \in \fa^l$. To show that all the $L_\Lambda$'s come  
from a globally defined $L$  
on all of $(AN)^l$, one uses the Moser Isotopy Lemma  
for Hamiltonian Poisson $K$--spaces proved in \cite{am:kashiwara}.  
What we present here is a collection of arguments from  
\cite{anton:jdg} \cite{amw:linear} \cite{KMT} and \cite{am:kashiwara}.

\subsection{Gauge transformation for Poisson structures} 
\label{sec_gauge} 
 
Recall that a Poisson structure on a manifold $M$ is a smooth section $\pi$ of  
$\wedge^2 TM$ such that $[\pi, \pi] = 0$, where $[ \, , \, ]$ is the Schouten 
bracket on the space of multi-vector fields on $M$. For a  
smooth section $\pi$ of $\wedge^2 TM$,  we will use $\pi^{\#}$ to denote  
the bundle map $\pi^\#: T^*M \to TM: \pi^\#(\alpha) = \pi(\cdot, \alpha)$ for all cotangent 
vectors 
$\alpha$. Similarly,  
for a $2$-form $\gamma$ on $M$, we will set 
$\gamma^{\#}: TM \to T^*M: \gamma^{\#}(v)= \gamma(\cdot, v)$ for all tangent vectors $v$. 
Suppose now that $\pi$ is a Poisson structure on $M$ and that $\gamma$ is a closed $2$--form 
on $M$. If the bundle map 
$1 + \gamma^{\#} \pi^{\#}: T^*M \to T^*M$ 
is invertible, the section $\pi^\prime$ of $\wedge^2 TM$ given by 
\begin{equation} 
\label{eq_gauge} 
(\pi^\prime)^{\#} = \pi^{\#} (1 + \gamma^{\#} \pi^{\#})^{-1}: \, \,  
T^*M \lrw TM 
\end{equation} 
is then a Poisson structure on $M$.   
The Poisson structure $\pi^\prime$  will be called the 
{\it gauge transformation} of 
$\pi$ by the closed $2$-from $\gamma$, and we write 
$\pi^\prime = {{\mathcal G}}_{\gamma}(\pi)$.  
It is clear from (\ref{eq_gauge}) that $\pi$ and $\pi^\prime$ have the same  
symplectic leaves. If $S$ is a common symplectic leaf, then the symplectic $2$-forms 
$\omega^\prime$ and $\omega$ coming from $\pi^\prime$ and $\pi$ differ by $i_{S}^{*} \gamma$, 
where $i_S: S \to M$ is the inclusion map. See \cite{alan-s:3} for more detail.

\subsection{The Poisson Lie groups $(K, s\pi_K)$ and $(AN, \bullet_s, \pi_{AN,s})$} 
\label{sec_family} 
 
The  
group $AN$ carries a distinguished Poisson structure $\pi_{AN}$. Indeed,  
let $\la \, , \, \ra$ be the imaginary part of the  
Killing form  
of $\fg$ and identify $\fk$  
with $(\fa + \fn)^*$ via $\la \, , \, \ra$. 
For $x \in \fk$, denote by $\bar{x}$ the right invariant $1$-form on $AN$  
defined by $x$. 
Let $x_{AN}$ be the generator of the action of $\exp(tx)$ on $AN$ according to 
(\ref{eq_K-on-AN}). 
Then the unique section $\pi_{AN}$ of $\wedge^2 T(AN) $ such that  
\begin{equation} 
\label{eq_piAN} 
\pi_{AN}^{\#}(\bar{x}) = x_{AN}, \hs \hs \forall x \in \fk 
\end{equation} 
is a Poisson bi-vector field on $AN$. The Poisson structure $\pi_{AN}$ makes 
$AN$ into a {\it Poisson Lie group} in the sense that the group multiplication  
map $AN \times AN \to AN: (b_1, b_2) \mapsto b_1b_2$ is a Poisson map, where $AN \times AN$ 
is equipped with the product Poisson structure $\pi_{AN} \times \pi_{AN}$.  
We refer to \cite{lu-we:poi} and \cite{lu:mom} 
for details on Poisson Lie groups and Poisson Lie group actions and to  
\cite{lu-ratiu}) for details on $\pi_{AN}$. In particular, the dual Poisson Lie group 
of $(AN, \pi_{AN})$ is $K$ together with the   
Poisson 
structure $\pi_K$, explicitly given by $\pi_K = \Lambda_{0}^{r} - 
\Lambda_{0}^{l}$, where  
\[ 
\Lambda_0 = {1 \over 2} \sum_{\alpha \in \Delta^+} X_\alpha \wedge Y_\alpha 
\in \wedge^2 \fk 
\] 
with $X_\alpha, Y_\alpha \in \fk$ given in (\ref{eq_fk}) and 
$\Lambda_{0}^{r}$ 
and $\Lambda_{0}^{l}$ being respectively the right and left invariant  
bi-vector fields on $K$ determined by $\Lambda_0$. 
It follows from (\ref{eq_piAN}) that the symplectic leaves of $\pi_{AN}$  
are precisely the orbits of the $K$--action on $AN$ given in (\ref{eq_K-on-AN}). 
 
Let now $d$ be an involutive automorphism of the Dynkin diagram of $\fg$, and let 
$\tau_d$ be the quasi--split real form of $\fg$ given in (\ref{eq_tau-d}). 
Recall that $\tau_d$ 
leaves $AN$ invariant and defines a group isomorphism on $AN$. It is easy to check   
(see also Section 2.3 of \cite{amw:linear}) that 
$\tau_d: (AN, \pi_{AN}) \to (AN, \pi_{AN})$ is anti-Poisson, i.e.,  
${\tau_d}_*\pi_{AN} =  
-\pi_{AN}.$ 
Similarly, $\tau_d: K \to K$ is anti--Poisson for $\pi_K$. 
We will denote the restrictions of $\tau_d$ to $K$ and to $AN$ both by $\tau_d$, and 
we will refer to $(K, \pi_K, \tau_d)$ and $(AN, \pi_{AN}, \tau_d)$ as a 
{\it dual pair} of {\it Poisson Lie groups with anti--Poisson involutions}. In the context 
of Poisson Lie groups, the $K$--action on $AN$ given in (\ref{eq_K-on-AN}) and the 
$AN$--action on $K$ given in (\ref{eq_AN-on-K-1}) (with $AN$ replacing $A_0N_0$ and 
$K$ replacing $K_0$) are respectively called the {\it dressing actions}.

As is noticed in \cite{anton:jdg} and \cite{KMT}, we in fact have a smooth family of Poisson 
Lie groups $(AN, \bullet_s, \pi_{AN,s})$ for $s \in \R$.  
 Indeed, for $s \in \R-\{0\}$, let $F_s: \fp \to \fp$ 
be the diffeomorphism $x \mapsto sx$, and let 
$I_s = E \circ F_s \circ E^{-1}:  AN \to AN.$ 
Let $\pi_{AN, s}$ be the Poisson bi-vector field on $AN$ such that $(I_s)_*  
\pi_{AN, s} = s \pi_{AN}$, or, 
\begin{equation} 
\label{eq_piANt} 
\pi_{AN, s}(b)  = s (I_{s}^{-1})_* \left(\pi_{AN}(I_s(b)) \right), \hs \hs 
b \in AN. 
\end{equation} 
Define the group structure ${\bullet}_s: AN \times AN \to AN$ by 
\[ 
b_1 \bullet_s b_2 := I_{s}^{-1}(I_s(b_1)I_s(b_2)), \hs \hs b_1, b_2 \in AN. 
\] 
Then since $s \neq 0$, the map 
$(AN \times AN, \, \pi_{AN, s} \times \pi_{AN, s}) \to (AN, \pi_{AN, s}): 
(b_1, b_2) \mapsto b_1b_2$ is a  
Poisson map, so $(AN, {\bullet}_s, \pi_{AN, s})$ is a Poisson Lie group 
for each $s \in \R$, $s \neq 0$. On the other hand,  
$\fp \cong \fk^*$ has the linear Poisson structure  
$\pi_{{\mathfrak p},0}$ defined by the Lie algebra $\fk$. 
Let ${\bullet}_0$ and $\pi_{AN, 0}$ be the pullbacks to $AN$ by 
$E^{-1}: AN \to \fp$ of the abelian 
group structure on $\fp$ and the Poisson structure $\pi_{{\mathfrak p},0}$ on $\fp$. 
Then we get a smooth family of Poisson Lie group structures $({\bullet}_s, \pi_{AN, s})$ 
on $AN$ for every $s \in \R$ (see \cite{anton:jdg} and \cite{KMT}).  
The dual Poisson Lie group of 
$(AN, {\bullet}_s, \piANt)$ is again the Lie group $K$ (with
 group structure independent on $s$) 
 but with the Poisson structure $s\pi_{K}$, if we identify again  
$\fk \cong (\fa + \fn)^*$ via 
the imaginary part of the Killing form of $\fg$. 
It is also clear that  
$\tau_d$ is a group isomorphism for 
$\bullet_s$ and is anti--Poisson for $\pi_{AN,s}$. Thus we get a dual pair 
of Poisson Lie groups $(K, s\pi_K, \tau_d)$ and $(AN, \bullet_s, \pi_{AN,s}, \tau_d)$ with  
anti--Poisson involutions for each $s \in \R$. 
 
\subsection{The Poisson structure $\pi_s$ on $(AN)^l$} 
 
As is noted in \cite{anton:jdg}, for each $s \in \R$, the Poisson structure $\pi_{AN,s}$ 
on $AN$ is related to $\pi_{AN}$ by a gauge transformation. Recall that for $x \in \fk  
\cong (\fa+\fn)^*$, $\bar{x}$ is the right invariant $1$--form on $AN$ defined by $x$. 
Let $l_x$ be the differential of the linear function $\xi \mapsto \la x, \xi \ra$ 
on $\fp$, and let $x_{{\mathfrak p}}$ be the vector field on $\fp$ generating the 
adjoint action of $\exp (tx) \in K$ on $\fp$. 
By Proposition 3.1 of \cite{amw:linear},  
there is a $1$-form $\beta$ on $\fp$ such that $\beta(0)=0$ and 
\[ 
-(d\beta)^{\#}(x_{{\mathfrak p}}) = E^* \bar{x} - l_x, \hs \hs \forall 
x \in \fk, 
\] 
Moreover, $(\tau_d)^* \beta = -\beta$ for every quasi-split real form  
$\tau_d$ given in (\ref{eq_tau-d}). 
For $s \in \R-\{0\}$, let $\beta_s ={1 \over s} (F_{s}^{*} \beta) $.  
Since $\beta(0) = 0$, the family $\beta_s$ extends smoothly to $\beta_0 = 0$.  
Let 
\[ 
\alpha = (E^{-1})^* \beta, \hs {\rm and} \hs 
\alpha_s = (E^{-1})^* \beta_s = {1 \over s} (I_{s}^{*}\alpha) 
\] 
for all $s \in \R$. Then it is easy to show that, for every $s \in \R$, 
\[ 
\pi_{AN, s} = {{\mathcal G}}_{d (\alpha - \alpha_s)}(\pi_{AN}) =  
{{\mathcal G}}_{-d\alpha_s}(\pi_{AN, 0}). 
\]  
 
Assume now that $l \geq 1$ is an integer. For each $s \in \R$, set  
\[ 
{\bf m}_s: \, (AN)^l \lrw AN: \, \, (b_1, b_2, \cdots, b_l) \Map 
b_1 \bullet_s b_2 \bullet_s \cdots \bullet_s b_l. 
\] 
By generalizing the ``linearization" procedure of Hamiltonian symplectic 
$(K, s\pi_K)$--spaces described in \cite{amw:linear} to the case of Poisson manifolds, 
one can show that  
\begin{equation} 
\label{eq_pis} 
\pi_s\, := \, {{\mathcal G}}_{d {\bf m}_{s}^{*} \alpha_s} (\pi_{AN, s} \times \pi_{AN, s} 
\times \cdots \times \pi_{AN, s}) 
\end{equation} 
is a well-defined Poisson structure on $(AN)^l$ for each $s \in \R$. Define the twisted diagonal action 
${{\mathcal T}}_{s}$ of $K$ on $(AN)^l$ by 
\begin{equation} 
\label{eq_twisted-s} 
k \Map {{\mathcal T}}_{s, k}  \, := \, \nu_{s}^{-1} \circ \delta_k \circ \nu_s, 
\end{equation} 
where, again, $\delta_k$ denotes the diagonal action of $k$ on $(AN)^l$ for the 
action of $K$ on $AN$ given in (\ref{eq_K-on-AN}), and $\nu_s \in {\rm Diff}((AN)^l)$ 
is given by 
\[ 
\nu_s(b_1, b_2, \cdots, b_l) = (b_1, \, \, b_1 \bullet_s b_2, \, \, \cdots, \, \,  
b_1 \bullet_s b_2 \bullet_s \cdots \bullet_s b_l). 
\] 
Note that ${{\mathcal T}}_{s, k}$ is ${{\mathcal T}}_{k}$ when $s=1$ and
is $\delta_k$ when $s=0$.
Then again it follows from \cite{amw:linear} that for each $s \in \R$, 
the action ${{\mathcal T}}_s$ of $K$ on $(AN)^l$ is Hamiltonian with respect to the 
Poisson structure $\pi_s$ with the map $E^{-1} \circ {\bf m}_s: (AN)^l \to \fp \cong \fk^*$ 
as a moment map. Moreover, for every quasi--split real form $\tau_d$ defined in (\ref{eq_tau-d}), 
the Cartesian product $(\tau_d)^l = \tau_d \times \tau_d \times \cdots 
\times \tau_d$ is anti--Poisson for $\pi_s$ for every $s \in \R$. 
 
\subsection{The Moser Isotopy Lemma} 
 
Let $U$ be a connected Lie group with 
Lie algebra $\fu$. Suppose that $\sigma_U$ is an involutive automorphism  of $U$  
with the corresponding  
involution on $\fu$ denoted by $\sigma_{{\mathfrak u}}$.  
Define $\sigma_{{\mathfrak u}^*} = -(\sigma_{{\mathfrak u}})^*$.  
If $(M, \pi_M, \Phi)$ is  
a Hamiltonian Poisson $U$--space, an anti-Poisson involution $\sigma_M$ of 
$(M, \pi_M)$ is said to be compatible with $\sigma_U$ if  
$\Phi \circ \sigma_M  = \sigma_{{\mathfrak u}^*} \circ \Phi.$ 
The following Moser Isotopy Lemma for Hamiltonian Poisson $U$--spaces with anti--Poisson 
involutions is proved in \cite{am:kashiwara}. 
See \cite{amw:linear} for the symplectic case.

\begin{prop} 
\label{prop_moser} Let $U$ be a connected compact semi-simple Lie group 
with Lie algebra $\fu$, and let  
$(M, \pi_s,  \Phi_s)$, $s \in \R$, be a smooth 
family of Hamiltonian Poisson $U$--spaces. Suppose that there  
exists a smooth family of $1$-forms $\epsilon_s$ on $M$ with $\epsilon_0 = 0$ 
such that $\pi_s = {{\mathcal G}}_{d \epsilon_s} \pi_0$ for 
every $s \in \R$.  
Assume also that $\pi_0$ has compact symplectic leaves. Then $(M, \pi_s, \Phi_s)$  
is 
isomorphic to $(M, \pi_0, \Phi_s)$ for every $s \in \R$ as  
Hamiltonian Poisson $U$--spaces, 
i.e., there exists 
$\psi_s \in {\rm Diff}(M)$ for $s \in \R$ with $\psi_0 = {\rm id}$, 
% CHANGED phi_0 IN ABOVE LINE TO  psi_0
such that for every $s \in \R$, 
 
1) $ \pi_s = {\psi_s}_*\pi_0$; 
%GOT RID OF INVERSE  
2) $\Phi_s \circ \psi_s = \Phi_0$. 
 
\noindent 
If $\sigma_U$ is an involutive automorphism on $U$, and if for each $s \in \R$, 
$\sigma_{M,s}$ is 
an anti-Poisson involution for $\pi_s$ compactible with $\sigma_U$ and 
such that $\sigma_{M, s}^{*}\dot{\epsilon}_s 
=-\dot{\epsilon}_s$,  
then $\psi_s$ can be chosen such that $\psi_s \circ \sigma_{M, 0} = \sigma_{M, s} 
\circ \psi_s$ for all $s \in \R$. 
\end{prop} 
 
\subsection{Proof of Theorem \ref{thm_amw}} 
 
Consider the Hamiltonian Poisson $K$--space $(M=(AN)^l, \pi_s, \Phi_s)$ with  
$\Phi_s = E^{-1} \circ {\bf m}_s$. The action of $K$ on $(AN)^l$ induced by 
$(\pi_s, \Phi_s)$ is the twisted diagonal action ${{\mathcal T}}_s$ given in (\ref{eq_pis}). 
From the definition of $\pi_s$, we know that 
$\pi_s = {{\mathcal G}}_{d \epsilon_s} \pi_0$, where 
\[ 
\epsilon_s = {\bf m}_{s}^{*} \alpha_s - \sum_{j=1}^{l} p_{j}^{*} \alpha_s 
\] 
with $p_j: (AN)^l \to AN$ denoting the projection to the $j$'th factor.  
For every quasi-split real form  
$\tau_d$ given in (\ref{eq_tau-d}), since $\tau_d$ is a group isomorphism for 
$(AN, \bullet_s)$, we have $\tau_{d}^{*} \epsilon_s = -\epsilon_s$, and thus 
$\tau_{d}^{*}\dot{\epsilon}_s 
=-\dot{\epsilon}_s$ for every $s \in \R$. Let $\sigma_{M,s} = (\tau_d)^l$ and let 
$\psi_s \in {\rm Diff}((AN)^l)$ 
be as in Proposition \ref{prop_moser}. Then 
$L:=(E^{-1})^l \circ {\psi}_{1}^{-1}: (AN)^l \to (\fp)^l$ is the diffeomorphism
 in Theorem \ref{thm_amw}. 
% CHANGED \PHI_1 TO \PSI_1 AGAIN. ALSO ADDED AN INVERSE TO \PSI_1 SINCE IT
% SEEMS NEEDED.
Indeed, 
\[
E\circ a \circ L = E\circ a \circ (E^{-1})^l \circ {\psi}_{1}^{-1} = 
E\circ \Phi_0 \circ {\psi}_{1}^{-1} = E\circ \Phi_1 = m,
\]
where the second equality follows from the identity
\[
a\circ (E^{-1})^l = E^{-1}\circ m_0,
\]
which is a trivial consequence of the fact that $m_0$ is the pullback
of addition by the map $E$.

\qed 

Sam Evens, Department of Mathematics, The University of Notre Dame;

Jiang-Hua Lu, Department of Mathematics, the University of Hong Kong.

{\it email address:}  evens.1@@nd.edu and jhlu@@maths.hku.hk


\begin{thebibliography}{99} 
 
\bibitem[A-B-V]{abv:l} 
Adams, J., Barbasch, D., and Vogan, D., {\em The Langlands classification 
and irreducible characters for real reductive groups}, 
Birkhauser, 1992. 
 
%\bibitem[A-V]{a-v:l} 
%Adams, J., and Vogan, D., L-groups, projective representations, 
%and the Langlands classification, {\it Amer. J. Math.},  
%{\bf 114} (1) (1991), 45 - 138. 
 
\bibitem[Al]{anton:jdg} 
Alekseev, A., On Poisson actions of compact Lie groups on symplectic manifolds,  
{\em J. Diff. Geom.} {\bf 45} (1997), 241 - 256. 
 
\bibitem[Al-Me-W]{amw:linear} 
Alekseev, A., Meinrenken, E., and Woodward, C., Linearization of Poisson  
actions and singular values of matrix products, {\em Ann. Inst. Fourier}, Grenoble,  
{\bf 51} (6) (2001), 1691 - 1717. 
 
\bibitem[Al-Me]{am:kashiwara} 
Alekseev, A., Meinrenken, E., Poisson geometry and the Kashiwara--Vergne conjecture,  
math.RT/0209346. 
 
\bibitem[Ar]{araki:satake} 
Araki, S., On root systems and an infinitesimal classification of 
irreducible symmetric spaces, {\em J. Mathematics, Osaka City University}, 
{\bf 13} (1) (1962), 1-34. 
 
%\bibitem[B-S]{b-reyer:subgroup} 
%Berenstein, A., and Sjamaar, R., Coadjoint orbits, moment polytopes, and the  
%Hilbert-Mumford criterion, {\em J. Amer. math. Soc.}, {\bf 13} (2000), 433 - 466. 
  
%\bibitem[B-Ru]{radko:gauge} 
%Bursztyn, H., and Radko, O., Gauge equivalence of Dirac structures and  
%symplectic groupoids, math.SG/0202099. 
 
%\bibitem[E-L2]{%:harm} 
%S. Evens and J.-H. Lu. Poisson harmonic forms, Kostant harmonic forms,  
%and the $S^1$-equivariant cohomology of $K/T$.  {\it Adv. Math.}, {\bf 142}: 171-220,    
%1999.  
 
%\bibitem[E-L2]{e-l:l1} 
%Evens, S., and Lu, J.-H., On the variety of 
%Lagrangian subalgebras, {\it Ann. E.N.S.} {\bf 34} (2001), 631--668. 
 
%\bibitem[F-Ra]{f-ratiu} 
%Flaschka, H., and Ratiu, T., A convexity theorem for Poisson actions of compact Lie groups, 
%{\em Ann. Sci. Ecole Norm. Sup.} {\bf 29} (6) (1996), 787 - 809. 
 
\bibitem[Fo]{philip:loci} 
Foth, P.,  A note on Lagrangian loci of quotients, math.SG/0303322. 
 
\bibitem[Fu]{fulton:survey} 
Fulton, W., Eigenvalues, invariant factors, highest weights, and Schubert calculus, 
{\em Bull. Amer. Math. Soc.} (N.S.) {\bf 37}(3) (2000), 209 - 249. 
 
%\bibitem[G-W]{viktor-alan} 
%Ginzburg, V., and Weintein, A., Lie-Poisson structure on some Poisson Lie groups, {\em 
%J. Diff. Geom.} {\bf 5}(2) (1992), 445 - 453. 
 
\bibitem[Ka-Mi-L1]{KML:1} 
Kapovich, M., Leeb, B., and Millson, J., Polygons in symmetric spaces and buildings, Preprint,  
2002. 
 
\bibitem[Ka-Mi-L2]{KML:2} 
Kapovich, M., Leeb, B., and Millson, J., The generalized triangle inequalities in 
symmetric spaces and buildings with applications to algebra, Preprint, 2002. 
 
\bibitem[Ka-Mi-T]{KMT} 
Kapovich, M., Millson, J., and Treloar, T., The symplectic geometry of 
polygons in hyperbolic $3$-space, {\em Asian J. math.} (Kodaira's issue), 
{\bf 4}(1) (2000), 123 - 164. 
 
 
%\bibitem[Kl1]{klya:linear} 
%Klyachko, A., Stable bundles, representation theory, and Hermitian operators,  
%{\em Select Math.} (N.S.) {\bf 4} (3) (1998), 419 - 445. 
 
\bibitem[Kl]{klya:thompson} 
Klyachko, A., Random walks on symmetric spaces and inequalities for matrix spectra, {\em Linear 
Algebra and Appl.} {\bf 319} (1-3) (2000), 37 - 59. 
 
 
\bibitem[Kn1]{knapp:examples} 
Knapp, A., {\em Representation theory of semi-simple groups}, Princeton University Press, 1986. 
 
%\bibitem[Kn2]{knapp:green} 
%Knapp, A., {\em Lie groups beyond an introduction}, Progress in Mathematics, 
%Birkhauser, 1996. 
 
\bibitem[Ku-M]{kumar-millson} 
Kumar, S., and Millson, J., The generalized triangular inequalities for rank $3$ symmetric 
spaces of non-compact type, math.SG/0303264. 
 
%\bibitem[He1]{helgason:1} 
%Helgason, S., {\em Differential geometry, Lie groups, and symmetric spaces}, 
%Academic press, 1978. 
 
\bibitem[L-W]{lu-we:poi} 
Lu, J.-H.,  and Weinstein, A.,  Poisson Lie groups, dressing 
transformations, and Bruhat decompositions. {\em J. Diff. Geom.},  
{\bf 31} (1990), 501-526. 
 
\bibitem[L]{lu:mom} 
Lu, J.-H., Moment mappings and reductions of Poisson Lie groups, {\em Proc. Seminaire 
Sud-Rhodanien de Geometrie}, MSRI series, Springer-Verlag 1991, 209 - 226. 
 
\bibitem[L-Ra]{lu-ratiu} 
Lu, J.-H., and Ratiu, T., On the non-linear convexity theorem of Kostant, 
{\em J. of Amer. Math. Soc.} {\bf 4}(2) (1991), 349 - 361. 
 
\bibitem[O-S]{o-reyer:pair}
O'Shea, L., and Sjamaar, R., Moment maps and Riemannian symmetric pairs,
{\em Math. Ann.}, {\bf 317} (3) (2000), 415 - 457.


%\bibitem[Sm]{samelson:notes}
%Samelson, H., {\em Notes on Lie algebras}, Universitext, Springer-Verlag, 1990.

%\bibitem[R-Sp]{r-s:involution}
%Richardson, R. W. and Springer, T. A., The Bruhat order on symmetric varieties, 
%{\em Geometriae Dedicata} {\bf 35} (1990), 389 - 436.

%\bibitem[Sa]{satake:compactification}
%Satake, I., On representations and compactifications of
%symmetric Riemannian spaces, {\em Ann. Math.} {\bf 71} (1960), 77 - 110.

\bibitem[Sl]{schli:symmetric}
Schlichtkrull, H., {\em Hyperfunctions and harmonic analysis on 
symmetric spaces}, Birkhauser, 1984.

\bibitem[Se-W]{alan-s:3}
Severa, P., and Weinstein, A., Poisson geometry with a $3$-form background,
{\em Proceedings of the international workshop on non-commutative geometry
and string theory}, Keio University (2001). Available also at
math.SG/0107133.


\bibitem[Wa]{wallach:book}
Wallach, N., {\em Real reductive groups I}, Academic Press, 1988.

%\bibitem[X]{xu:involution}
%Xu, P., Dirac submanifolds and Poisson involutions, preprint.

\end{thebibliography}
\end{document}